\documentstyle[11pt]{article}

\title{NONSTANDARD GRAPHS, REVISED}
\author{A. H. Zemanian}
\date{}

\begin{document}
\newcommand{\N} {I \kern -4.5pt N}
\newcommand{\M} {I \kern -4.5pt M}
\newcommand{\R} {I \kern -4.5pt R}
\maketitle
\baselineskip21pt

{\ Abstract --- This is a revision of the paper archived 
previously on August 22, 2002.  It corrects a mistake in Sec. 8 
concerning eccentricities of graphs.  

From any given sequence of finite or infinite graphs, 
a nonstandard graph is constructed.
The procedure is similar to an ultrapower construction of an internal
set from a sequence of subsets of the real line, but now the individual
entities are the vertices of the graphs instead of real numbers. 
The transfer principle is then invoked to extend several graph-theoretic
results to the nonstandard case. After incidences and adjacencies between 
nonstandard vertices and edges are defined, 
several formulas regarding numbers of vertices and edges,
and nonstandard versions of Eulerian graphs, Hamiltonian
graphs, and a coloring theorem are established for these nonstandard graphs.

\vspace{.05in}
}

1. Introduction

In the book \cite{go} (se Sec. 19.1), R. Goldblatt 
obtained a nonstandard graph $\langle ^{*}\!V,^{*}\!E \rangle$
as an enlargement of a conventional graph $\langle V,E\rangle$
defined as a set $V$ of vertices and a set $E$ of edges 
determined by a symmetric irreflexive relationship of $V$.
The transfer principle was then used to establish
in a nonstandard way the standard theorem 
that, if every finite subgraph of a conventional infinite graph
$G$ has a coloring with finitely many colors, 
then $G$ itself has such a finite coloring.

On the other hand, in several recent works (see \cite{ze00a},
\cite{ze00b}, and the references therein), 
the idea of nonstandard transfinite graphs and networks
was introduced and investigated.  The basic idea in those works was to 
start with a given transfinite graph, to reduce it to a finite graph 
by shorting and opening edges,
and then to obtain an expanding 
sequence of finite graphs by restoring edges sequentially.  
If all this is done in an appropriate fashion, it may happen that the 
sequence of finite graphs fills out and restores the original 
transfinite graph once the restoration process is completed.  If in addition
there is an assignment of electrical parameters to the edges, 
the final result may be sequences of edge currents and edge
voltages, from which nonstandard 
hyperreal currents and voltages can be derived.
The latter hyperreal quantities
will then automatically satisfy Kirchhoff's laws, even though 
Kirchhoff's laws may on occasion be violated in the original 
transfinite network when standard real numbers are used---an 
important advantage of this nonstandard approach.

However, this is only a partial construction of a nonstandard graph in the
sense that the completion of the restoration 
process---if successful---results in the original standard transfinite graph.
The sequence of restorations only provides a means of constructing hyperreal 
currents and voltages satisfying Kirchhoff's laws.  Another approach 
might start with an arbitrary sequence of conventional 
(finite or infinite) graphs and 
construct from that a nonstandard graph in much the same way as
an internal set in the hyperreal line $^{*}\!\R$ 
is constructed from a given sequence of 
subsets of the real line $\R$, that is, by means of an 
ultrapower construction \cite[page 36]{go}.  In this case, the resulting
nonstandard graph has nonstandard edges and nonstandard
vertices, and it thereby is much 
different from those of the prior works \cite{ze00a}, \cite{ze00b}.
Also, the present approach differs from Goldblatt's result
in that the graphs of the sequence need not be vertex-induced 
subgraphs of a given graph.

Our objective in this work is to develop this latter approach to 
nonstandard graphs.  The individual elements are chosen to be the vertices 
of the graphs along with all the natural numbers.  These individuals
are not sets by assumption and therefore have no members.
All the other standard and nonstandard entities are derived 
from these individuals.  For example, an edge is defined to be a pair of 
vertices, which is one of the conventional ways of defining a graph.
Then, certain equivalence classes of sequences of standard vertices
are defined to be nonstandard vertices, and these then 
yield nonstandard edges as certain equivalence classes of sequences of 
standard edges.
In this approach, there are no multiedges (i.e., no parallel
edges) and no self-loops (i.e., no edge consisting of just a single 
vertex).

After setting up our nonstandard graphs using an ultrapower approach, 
we invoke the transfer principle to lift several standard
graph-theoretic results into a nonstandard setting.  For example,
the relationship between the number of vertices, the number of edges,
and the cyclomatic number of a standard finite connected 
graph continues to hold for nonstandard graphs except that these 
numbers are replaced by hypernatural numbers.  Similarly, standard 
theorems concerning 
Eulerian graphs, Hamiltonian graphs, and a coloring theorem are 
also extended to the nonstandard setting.  By virtue of the 
transfer principle, this only requires that the 
standard theorems be stated as 
sentences in symbolic logic, which are then transferred to 
appropriate sentences for nonstandard graphs.

In the following, $|A|$ will denote the cardinality of a set $A$.  Also,
$\N=\{0,1,2,\ldots\}$ is the set of all natural numbers,
and $\R$ is the set of real numbers.  Thus, $^{*}\!\N$ is the
set of hypernaturals, and $^{*}\!\R$ is the set of hyperreals.

2.Nonstandard Graphs

A standard graph $G$ is a conventional (finite or infinite) graph
$G=\{X,B\}$, where $X$ is the set 
of its vertices 
and $B$ is the set of its edges.  Each
edge $b\in B$ is a two-element set $b=\{ x,y\}$
with $x,y\in X$ and $x\neq y$;
$b$ and $x$ are said to be incident and so, too, are $b$
and $y$.  Also, $x$ and $y$ are said to be adjacent through $b$.  

Next, let $\langle G_{n}\!:n\in\N\rangle$ be a given
sequence of graphs. For each $n$, we have $G_{n}=\{X_{n},
B_{n}\}$, where $X_{n}$ is the set of edges 
and $B_{n}$
is the set of vertices.  We allow $X_{n}\cap X_{m}\neq\emptyset$
so that $G_{n}$ and $G_{m}$ may be subgraphs of a larger graph.
In fact, we may have $X_{n}=X_{m}$ and $B_{n}=
B_{m}$ for all $n,m\in\N$ so that $G_{n}$
may be the same graph for all $n\in\N$.
Furthermore, let ${\cal F}$ be a chosen 
nonprincipal ultrafilter on $\N$ \cite[pages 18-19]{go}.

In the following, $\langle x_{n}\rangle =\langle x_{n}\!: n\in\N\rangle$
will denote a sequence 
of vertices with $x_{n}\in X_{n}$ for all $n\in\N$.  
A nonstandard vertex $^{*}\!x$ is an equivalence class
of such sequences of vertices, where two such sequences 
$\langle x_{n}\rangle$ and $\langle y_{n}\rangle$ are taken
to be equivalent if $\{n\!: x_{n}=y_{n}\}\in{\cal F}$,
in which case we write ``$\langle x_{n}\rangle=\langle
y_{n}\rangle$ a.e.''
or say that $x_{n}=y_{n}$ ``for almost all $n$.''  We also
write $^{*}\!x=[x_{n}]$, where it is understood 
that the $x_{n}$ are the 
members of any one sequence in the equivalence class.

That this truly defines an equivalence class can be shown as follows.
Reflexivity and symmetry being obvious, consider transitivity:
Given that $\langle x_{n}\rangle=\langle y_{n}\rangle$ a.e.
and that $\langle y_{n}\rangle = \langle z_{n}\rangle$ a.e.,
we have $N_{xy}=\{n\!: x_{n}=y_{n}\}\in{\cal F}$ and 
$N_{yz}=\{n\!: y_{n}=z_{n}\}\in{\cal F}$.
By the properties of the ultrafilter, $N_{xy}\cap N_{yz}\in{\cal F}$.  
Moreover, $N_{xz}=\{n\!: x_{n}=z_{n}\}\supseteq
(N_{xy}\cap N_{yz})$.  Therefore, $N_{xz}\in{\cal F}$. 
Hence, $\langle x_{n}\rangle=\langle z_{n}\rangle$ a.e.; 
transitivity holds.  We let $^{*}\!X$ denote the set of 
nonstandard vertices.

Next, we define the nonstandard edges:  Let $^{*}\!x=[x_{n}]$ and 
$^{*}\!y=[y_{n}]$ be two nonstandard vertices.  This time, let 
$N_{xy}=\{n\!:\{x_{n},y_{n}\}\in B_{n}\}$ and 
$N^{c}_{xy}=\{n\!: \{x_{n},y_{n}\}\not\in B_{n}\}$.
Since ${\cal F}$ is an ultrafilter, exactly one of $N_{xy}$ and 
$N^{c}_{xy}$ is a member of ${\cal F}$.  If it is $N_{xy}$, then
$^{*}\! b=[\{x_{n},y_{n}\}]$ is 
defined to be a nonstandard edge;  that is, $^{*}\! b$ 
is an equivalence class of sequences $\langle b_{n}\rangle$
of edges $b_{n}=\{x_{n},y_{n}\}\in B_{n}$,
$n=0,1,2,\ldots\;$.  In this case, we also write $^{*}\!x,^{*}\!y
\in\,^{*}\!b$ and $^{*}\! b =\{ ^{*}\!x,\,^{*}\!y\}$.
We let $^{*}\!B$ denote
the set of nonstandard edges.  On the other hand, 
if $N^{c}_{xy}\in{\cal F}$, then $[\{ x_{n},y_{n}\}]$ is 
not a nonstandard edge.  

We shall now show that this definition is independent of the 
representatives chosen for the vertices.  Let $[\{ x_{n},y_{n}\}]$
and $[\{v_{n}, w_{n}\}]$ represent the same 
nonstandard edge.
We want to show that, if $\langle x_{n}\rangle=\langle 
v_{n}\rangle$ a.e., then $\langle y_{n}\rangle=
\langle w_{n}\rangle$ a.e.  Suppose $\langle y_{n}\rangle
\neq\langle w_{n}\rangle$ a.e. Then, $\{n\!: x_{n}=v_{n}\}
\cap\{ n\!:y_{n}\neq w_{n}\}\in {\cal F}$.
Thus, there is at least one $n$ for which the three vertices $x_{n}=
v_{n}$, $y_{n}$, and $w_{n}$ are all incident to the same 
standard edge---in violation of the definition of a edge.
Similarly, if all of $\langle x_{n}\rangle$, $\langle y_{n}\rangle$,
$\langle v_{n}\rangle$, $\langle w_{n}\rangle$ are different a.e.,
then there would be a standard edge having four incident vertices---again
a violation.

Next,  we show that we truly have an equivalence relationship
for the set of all sequences of standard edges.
Reflexivity and symmetry being obvious again, consider transitivity:
Let $^{*}\!b=[\{ x_{n},y_{n}\}]$, 
$^{*}\!\tilde{b}=[\{\tilde{x}_{n},\tilde{y}_{n}\}]$,
$^{*}\!\acute{b}=[\{\acute{x}_{n},\acute{y}_{n}\}]$,
and assume that $^{*}\!b=\, ^{*}\!\tilde{b}$ and 
$^{*}\!\tilde{b}=\, ^{*}\!\acute{b}$.  We want to show that
$^{*}\!b=\, ^{*}\!\acute{b}$.  We have $N_{b \tilde{b}}=\{n\!:
\{ x_{n}, y_{n}\}=\{\tilde{x}_{n},
\tilde{y}_{n}\}\}\in {\cal F}$ and
$N_{\tilde{b}\acute{b}}=\{ n\!: \{\tilde{x}_{n},
\tilde{y}_{n}\}=\{\acute{x}_{n}, \acute{y}_{n}\}\}\in {\cal F}$.
Moreover, $N_{b\acute{b}}=\{n\!:\{x_{n},y_{n}\}=
\{\acute{x}_{n},\acute{y}_{n}\}\}\supseteq
(N_{b\tilde{b}}\cap N_{\tilde{b}\acute{b}})\in{\cal F}$.
Therefore, $N_{b\acute{b}}\in{\cal F}$. Thus, $^{*}\!b=\,
^{*}\!\acute{b}$, as desired.

Finally, we define a nonstandard graph $^{*}\!G$
to be the pair $^{*}\!G=\{^{*}\!X,\,^{*}\!B\}$.

Let us now take note of a special case that arises when all the 
$G_{n}$ are the same standard graph $G=\{X,B\}$.  In this case, we call 
$^{*}\!G=\{^{*}\!X,\,^{*}\!B\}$ an {\em enlargement} of $G$,
in conformity with an ``enlargement'' $^{*}\!A$ of a subset $A$ of 
$\R$ \cite[pages 28-29]{go}.  If in addition $G$ is a finite graph, 
each vertex 
$x\in\,^{*}\!X$ can be identified with a vertex of $X$
because the enlargement of a finite set equals the set itself.
Similarly, every branch $b\in\,^{*}\!B$ can be identified with a branch in 
$B$.  In short, $^{*}\!G=G$.  On the other hand, if $G$ 
is a conventionally infinite graph, $X$ is an infinite set and its 
enlargement $^{*}\!X$ has more elements, namely, nonstandard vertices
that are not equal to standard vertices, (i.e., 
$^{*}\!X\setminus X$ is not empty).  Similarly, $^{*}\!B\setminus B$
is not empty too.  In short, $^{*}\!G$ is a proper enlargement of $G$.

{\bf Example 2.1.}  Let $G$ be a one-way infinite path $P$:
\[ P\;=\;\langle x_{0},b_{0},x_{1},b_{1},x_{2},b_{2},\ldots\rangle \]
Also, let $G_{n}=G$ for all $n\in\N$, and let $^{*}\!G=[G_{n}]
=\{^{*}\!X,\,^{*}\!B\}$.  Next, let $\langle k_{n}\!: n\in\N\rangle$
be any sequence of natural numbers, and set $^{*}\!x=[x_{k_{n}}]$ and 
$^{*}\! y=[x_{k_{n}+1}]$.  Then, $^{*}\!x$ and $^{*}\!y$
are two vertices in the enlargement $^{*}\!G$ of $G$, and 
$^{*}\!b=\{^{*}\!x,^{*}\!y\}=[\{x_{k_{n}},x_{k_{n}+1}\}]$
is an edge in $^{*}\!G$.  On the other hand, if 
$\langle m_{n}\!: n\in\N\rangle$ is another sequence of natural
numbers with $m_{n}>1$, then $^{*}\!z=[x_{k_{n}+m_{n}}]$
is another nonstandard vertex in $^{*}\!G$ different from
$^{*}\!x$ and $^{*}\!y$ and appearing after $^{*}\!y$ in
the enlarged path.  Moreover, $[\{x_{k_{n}},x_{k_{n}+m_{n}}\}]$
is not a nonstandard edge.  In this way, no vertex or edge repeats in 
the enlarged path.  $\Box$

Another special case arises when almost all the $G_{n}$ are 
(possibly different) finite graphs.  Again in conformity with the 
terminology used for hyperfinite internal subsets of $^{*}\!R$
\cite[page 149]{go}, we will refer to the resulting nonstandard graph
$^{*}\!G$ as a hyperfinite graph.\footnote{This should not be confused 
with a hypergraph---an entirely different object \cite{ber}.}
As a result, we can lift many theorems concerning finite graphs
to hyperfinite graphs.  It is just a matter of writing the standard theorem 
in an appropriate form using symbolic logic and then applying
the transfer principle.  We let $^{*}\!G_{f}$ denote the set of hyperfinite 
graphs.

3. Incidences and Adjacencies between Vertices and Edges

Let us now define these ideas for nonstandard graphs both in terms of 
an ultrapower construction 
and by transfer of appropriate symbolic sentences.
In the subsequent sections, we will usually confine ourselves to the
transfer principle.  We henceforth drop the asterisks when denoting 
nonstandard vertices and edges.  These are specified as members of 
$^{*}\! X$ and $^{*}\! B$ respectively.

{\em Incidence between a vertex and an edge:}  Given a sequence 
$\langle G_{n}\!: n\in\N\rangle$ of standard graphs $G_{n}=
\{X_{n},B_{n}\}$, a nonstandard vertex $x=[x_{n}]\in\,^{*}\! X$ 
and a nonstandard edge $b=[b_{n}]\in\,^{*}\! B$ are said to be
{\em incident} if $\{n\!: x_{n}\in b_{n}\}\in{\cal F}$, where as 
always the nonprincipal ultrafilter $\cal F$ is understood to be 
chosen and fixed.  

On the other hand, we can define incidence between a standard vertex $x\in X$
and a standard edge $b\in B$ for the graph $G=\{X,B\}$ through the 
symbolic sentence
\[ (\exists\; x\in X)\;(\exists\; b\in B)\;(x\in b) \]
By transfer, we have that $x\in\,^{*}\! X$ and $b \in\,^{*}\! B$
are incident when the following sentence is true.
\[ (\exists\; x\in\,^{*}\! X)\;(\exists\; b\in\,^{*}\! B)\;(x\in b) \]
These are equivalent definitions.

{\em Adjacency between vertices:}  For a standard graph $G=\{X,B\}$,
two vertices $x,y\in X$ are called {\em adjacent} and we write 
$x\diamond y$ if the following sentence on the right-hand side of 
$\leftrightarrow$ is true.
\[ x\diamond y\;\leftrightarrow\;(\exists\; x,y\in X)\;(\exists\; b\in B)\;(b=\{x,y\}) \]
By transfer, this becomes for a nonstandard graph
$^{*}\! G=\{^{*}\! X,\,^{*}\! B\}$
\[ x\diamond y\;\leftrightarrow\;(\exists\; x,y\in\,^{*}\! X)\;(\exists\; b\in\,^{*}\! B)\;(b=\{x,y\}) \]
Alternatively, under an ultrapower construction, we have
for $^{*}\! G=[G_{n}]=[\{X_{n},B_{n}\}]$ that 
$x=[x_{n}]\in\,^{*}\! X$ and $y=[y_{n}]\in\,^{*}\! X$ are adjacent
(i.e., $x\diamond y$) if there exists a $b=[b_{n}]\in\,^{*}\! B$ 
such that $\{n\!: b_{n}=\{x_{n},y_{n}\}\}\;\in {\cal F}$.

{\em Adjacency between edges:}  For a standard graph, two edges 
$b,c\in B$ are called {\em adjacent} and we write $b\bowtie c$
when the following sentence on the right-hand side of $\leftrightarrow$
is true.
\[ b\bowtie c\;\leftrightarrow\;(\exists\; b,c\in B)\;(\exists\; x\in X)\;(x\in b \wedge x\in c) \]
By transfer, we have for nonstandard edges $b$ and $c$
\[ b\bowtie c\;\leftrightarrow\;(\exists\; b,c\in\,^{*}\! B)\;(\exists\; x\in\,^{*}\! X)\;(x\in b \wedge x\in c) \]
Under an ultrapower approach, we would have $b=[b_{n}]$ and 
$c=[c_{n}]$ are adjacent nonstandard edges when there exists a nonstandard 
vertex $x=[x_{n}]$ such that 
\[ \{n\!: x_{n}\in b_{n}\; {\rm and}\;
x_{n}\in c_{n}\}\;\in\;{\cal F}. \]

4. Nonstandard Hyperfinite Paths and Loops

Again, we start with a standard graph $G=\{ X,B\}$. Remember that since $B$ 
is a set of two-element subsets of $X$, all edges are distinct 
(i.e., there are no multiedges) and there are no self-loops.
A {\em finite path} $P$ in $G$ is defined by the sentence
\[ (\exists\; k \in \N\setminus\{0\})\;(\exists\; x_{0},x_{1},\ldots,x_{k}\in X)\;
(\exists\; b_{0},b_{1},\ldots b_{k-1}\in B) \]
\begin{equation}
(x_{0}\in b_{0}\,\wedge\, b_{0}\ni x_{1}\,\wedge\, x_{1}\in b_{1}\,\wedge\,
b_{1}\ni x_{2}\,\wedge\,\ldots\,\wedge\, x_{k-1}\in b_{k-1}\,\wedge\, 
b_{k-1}\ni x_{k})  \label{4.1}
\end{equation}
That all the vertices and edges in $P$ are distinct is implied by the
fact that those $k$ vertices in $X$ and those $k-1$ edges in 
$B$ are perforce all distinct.
The {\em length} $|P|$ of $P$ is the number of edges
in $P$; thus, $|P|=k$.

We may apply transfer to (\ref{4.1}) to get the following definition 
of a {\em nonstandard path} $ ^{*}\! P$.
\[ (\exists\; k \in\,^{*}\! \N\setminus\{0\})\;(\exists\; x_{0},x_{1},\ldots,x_{k}\in\,^{*}\! X)\;
(\exists\; b_{0},b_{1},\ldots b_{k-1}\in\,^{*}\! B) \]
\begin{equation}
(x_{0}\in b_{0}\,\wedge\, b_{0}\ni x_{1}\,\wedge\, x_{1}\in b_{1}\,\wedge\,
b_{1}\ni x_{2}\,\wedge\,\ldots\,\wedge\, x_{k-1}\in b_{k-1}\,\wedge\, 
b_{k-1}\ni x_{k})  \label{4.2}
\end{equation}
In this case, the {\em length} $| ^{*}\! P|$ equals $k\in\,^{*}\! \N\setminus\{0\}$;
in general, $k$ is now a positive hypernatural number.  We therefore
call $^{*}\! P$ a {\em hyperfinite path}.  To view this fact in terms of an
ultrapower construction of $^{*}\! G=[G_{n}]$, note that the 
$G_{n}$ may be finite graphs growing in size or indeed be conventionally 
infinite graphs.  Thus, $^{*}\! P$ may have 
an unlimited length, that is, its length may be a member of 
$^{*}\! \N\setminus \N$. 

A standard loop is defined as is a standard path except that the 
first and last vertices are required to be the same.  Upon
applying transfer, we get the following definition of a {\em nonstandard loop}.
\[ (\exists\; k\in\,^{*}\! \N\setminus\{0\})\;(\exists\; x_{0},x_{1},\ldots,x_{k-1}\in\,^{*}\! X)\;
(\exists\; b_{0},b_{1},\ldots b_{k-1}\in\,^{*}\! B) \]
\begin{equation}
(x_{0}\in b_{0}\,\wedge\, b_{0}\ni x_{1}\,\wedge\, x_{1}\in b_{1}\,\wedge\,
b_{1}\ni x_{2}\,\wedge\,\ldots\,\wedge\, x_{k-1}\in b_{k-1}\,\wedge\, 
b_{k-1}\ni x_{0})  \label{4.3}
\end{equation}

5. Connected Nonstandard Graphs

A standard graph $G=\{ X,B\}$ is called {\em connected} if, for every two
vertices $x$ and $y$ in $G$, there is a finite path (\ref{4.1})
terminating at those vertices, that is, $x_{0}=x$ and $x_{k}=y$.
Let $\cal C$ denote the set of connected standard graphs, and let
${\cal P}(G)$ be the set of all finite paths in a given standard 
graph $G$.  Also, for any $P\in{\cal P}(G)$, let $x_{0}(P)$
and $x_{k}(P)$ denote the first and last nodes of $P$ in accordance
with (\ref{4.1});  $k$ depends upon $P$.
Then, the connectedness of $G$ is defined 
symbolically by the truth of the following sentence to the right of 
$\leftrightarrow$.
\begin{equation}
G\in {\cal C}\; \leftrightarrow\;(\forall\; x,y\in X)\;(\exists\; P\in{\cal P}(G))\; 
((x_{0}(P)=x)\,\wedge\, (x_{k}(P)=y))  \label{5.1}
\end{equation}
By transfer, we obtain the definition of the set $ ^{*}\! {\cal C}$
of all {\em connected} nonstandard graphs:  For $ ^{*}\! G=
\{ ^{*}\! X, ^{*}\! B\}$, for $ ^{*}\! {\cal P} ( ^{*}\! G)$
being the set of nonstandard paths $ ^{*}\! P\in\,^{*}\! G$,  
and for $x_{0}( ^{*}\! P)$ and $x_{k}( ^{*}\! P)$ being the first
and last vertices of $ ^{*}\! P$, we have
\begin{equation}
^{*}\! G\in\,^{*}\! {\cal C} \leftrightarrow\;
(\forall\; x,y\in\,^{*}\! X)\;(\exists\; ^{*}\! P\in\,^{*}\! {\cal P}( ^{*}\! G))\;
((x_{0}( ^{*}\! P)=x)\,\wedge\, (x_{k}( ^{*}\! P)=y))  \label{5.2}
\end{equation}
Here, $k\in\,^{*}\! \N\setminus\{0\}$ as in (\ref{4.2}).

Let us explicate this still further in terms of an ultrapower construction
of $ ^{*}\! G =[G_{n}]$ from an equivalence class of sequences of 
(possibly infinite) graphs, $\langle G_{n}\rangle$ being one of those 
sequences.  With $ ^{*}\! P=[P_{n}]$ denoting a nonstandard path 
obtained similarly from a representative sequence 
$\langle P_{n}\rangle$ of finite paths, $P_{n}$ being in $G_{n}$, 
we let $x_{0}(P_{n})$ and $x_{k}(P_{n})$ denote the first and last vertices
of $P_{n}$. (Thus, $k$ also depends on $n$ of course.)  Then, 
$x_{0}(^{*}\! P)=[x_{0}(P_{n})]$ and $x_{^{*}\! k}(^{*}\!P)=
[x_{k}(P_{n})]$ are the first and last nonstandard vertices of 
$^{*}\! P$. (In (\ref {5.2}), $^{*}\! k$ is denoted by 
$k\in\,^{*}\! \N\setminus\{0\}$.) Then, $^{*}\! G$ is called {\em connected} 
if and only if, given any nonstandard vertices $^{*}\! x=[x_{n}]$
and $^{*}\! y=[y_{n}]$ in $^{*}\! G$, we have that, for almost all $n$, 
there exists a finite path 
$P_{n}$ terminating at $x_{n}$ and $y_{n}$.  This can be restated by 
saying that there exists a hyperfinite path $^{*}\! P$ in 
$^{*}\! G$ terminating at $^{*}\! x$ and $^{*}\! y$.

Later on, we will need a special case of $^{*}\! {\cal C}$:
Let ${\cal C}_{f}$ denote the subset of ${\cal C}$ consisting of all 
finite connected standard graphs $G=\{X,B\}$,
where $|X|\in\N\setminus\{0,1\}$, 
$|B|\in \N\setminus \{0\}$.  Then, $^{*}\! {\cal C}_{f}$
is the subset of $\cal C$ obtained by lifting ${\cal C}_{f}$
through transfer to a subset of $^{*}\! {\cal C}$.  In this case,
\begin{equation}
^{*}\! G\in\,^{*}\! {\cal C}_{f}\; \leftrightarrow\;(^{*}\! G=
\{^{*}\! X,\,^{*}\! B\}\in\,^{*}\!{\cal C})\,\wedge\,
(|^{*}\! X|\in\,^{*}\!\N\setminus\{0,1\}\,\wedge\,|^{*}\! B|\in\,
^{*}\! \N\setminus\{0\}).  \label{5.3}
\end{equation}
We call such a $ ^{*}\! G$ a nonstandard {\em hyperfinite} connected graph.

6. Nonstandard Subgraphs

If $A$ and $C$ are sets of vertices with $A\subseteq C$
, we get upon transfer 
the following definition for sets of nonstandard vertices.
\[ ^{*}\! A\subseteq\, ^{*}\! C\; \leftrightarrow \;
(\forall\; x\in\,^{*}\! A)\;(x\in\,^{*}\! C) \]

Our purpose in this section is to define a nonstandard subgraph 
$ ^{*}\! G_{s}$ of a given nonstandard graph $ ^{*}\! G$.
We let $\cal G$ denote the set of all standard graphs.
By definition, $G_{s}=\{X_{s},B_{s}\}$ is a (vertex induced)
subgraph of $G=\{X,B\}\in{\cal G}$ if $X_{s}\subseteq X$ and $B_{s}$ 
is the set of those edges in $B$
that are each incident to two vertices in $X_{s}$.
Let us denote the set of all such subgraphs of $G$ by ${\cal G}_{s}(G)$.
Then, symbolically $G_{s}$ is defined by 
\[ G_{s}\in{\cal G}_{s}(G)\; \leftrightarrow\; \]
\[ (\exists\; G_{s}=\{X_{s},B_{s}\}\in{\cal G})\;
(\exists\; G=\{X,B\}\in{\cal G}) \]
\[ ((X_{s}\subseteq X)\,\wedge\,(B_{s}=\{b=\{x,y\}\in B\!: x,y\in X_{s}\})). \]
By transfer, we get the definition of a {\em nonstandard subgraph}
$ ^{*}\! G_{s}$ of a given nonstandard graph $ ^{*}\! G$.
In this case, $ ^{*}\! {\cal G}$ denotes the set of all 
nonstandard graphs, and $ ^{*}\! {\cal G}_{s}( ^{*}\! G)$ denotes
the set of all nonstandard subgraphs of a given 
$ ^{*}\! G\in  ^{*}\! {\cal G}$.
\[  ^{*}\! G_{s}\in\,^{*}\!{\cal G}_{s}( ^{*}\! G)\; \leftrightarrow\; \]
\[ (\exists\;  ^{*}\! G_{s}=\{ ^{*}\! X_{s},\, ^{*}\! B_{s}\}\in\,^{*}\! {\cal G})\;
(\exists\; ^{*}\! G=\{ ^{*}\! X, ^{*}\! B\}\in\,^{*}\! {\cal G}) \]
\[ (( ^{*}\! X_{s}\subseteq\,^{*}\! X)\,\wedge\,( ^{*}\! B_{s}=\{ b=\{x,y\}
\in\,^{*}\! B\!: x,y\in\,^{*}\! X_{s}\}) \]

7. Nonstandard Trees

The symbols ${\cal G}$, ${\cal G}_{s}(G)$,
${\cal C}$, ${\cal C}_{f}$, and their nonstandard
counterparts have been defined above.  Now, let ${\cal L}(G)$
be the set of all loops in the standard graph $G$, and let $\cal T$ 
be the set of all standard trees. Then, a tree $T=\{ X_{T},B_{T}\}$ 
can be defined symbolically by 
\begin{equation}
T\in {\cal T}\; \leftrightarrow\; (\exists\; T\in{\cal C})\;
(\neg(\exists\; L\in{\cal L}(T)))  \label{7.1}
\end{equation}
To transfer this, we let $ ^{*}\! {\cal L}( ^{*}\! G)$ be the set of all
nonstandard loops in a given nonstandard graph $ ^{*}\! G$, as defined
by (\ref{4.3}), and we let $ ^{*}\! {\cal T}$ denote the set
of {\em nonstandard trees} $ ^{*}\! T$, defined as follows:
\begin{equation}
 ^{*}\! T\in\,^{*}\! {\cal T}\; \leftrightarrow\;
(\exists\; ^{*}\! T\in\,^{*}\! {\cal C})\;
(\neg(\exists\;  ^{*}\! L\in\,^{*}\! {\cal L}( ^{*}\! T)))  \label{7.2}
\end{equation}

Next, let ${\cal T}_{sp}(G)$ be the set of all spanning trees in 
a given finite connected standard graph $G=\{X,B\}$.
That $T$ is such a spanning tree can be expressed symbolically as follows.
\[ T\in{\cal T}_{sp}(G)\; \leftrightarrow\; \]
\begin{equation}
(\exists\; G=\{X,B\}\in{\cal C}_{f})\;
(\exists\; T=\{X_{T},B_{T}\}\in {\cal T})\;
((T\in{\cal G}_{s}(G))\;\wedge\;(|X|=|X_{T}|))  \label{7.3}
\end{equation}
By transfer, we have the set $ ^{*}\! {\cal T}_{sp}( ^{*}\! G)$
of all {\em spanning trees} of a given hyperfinite connected nonstandard 
graph $ ^{*}\! G$, defined as follows:
\[  ^{*}\! T\in\,^{*}\! {\cal T}_{sp}( ^{*}\! G)\; \leftrightarrow\; \]
\begin{equation}
(\exists\; ^{*}\! G=\{ ^{*}\!X,\,^{*}\!B\}\in\,^{*}\! {\cal C}_{f})\;
(\exists\; ^{*}\! T=\{ ^{*}\! X_{T},\, ^{*}\! B_{T}\}\in\,^{*}\! {\cal T})\;
(( ^{*}\! T
\in\,^{*}\! {\cal G}_{s}( ^{*}\! G))\;\wedge\;
(| ^{*}\! X|=| ^{*}\! X_{T}|))  \label{7.4}
\end{equation}

8. Some Numerical Formulas

With these symbolic definition in hand, we can now lift some 
standard formulas regarding numbers of vertices and edges
into a nonstandard setting.  For example, if $p$, $q$, and $r$
are the number of vertices, the number of edges, and the 
cyclomatic number respectively of a given connected finite graph,
then $r=q-p+1$.  Symbolically, this can be stated as follows.  Again,
we use the notation $G=\{X,B\}$ for a graph and 
$T=\{X_{T},B_{T}\}$ for a tree.
\[(\forall\; p,q,r\in\N)\;(\forall\; G\in{\cal C}_{f})\;(\forall\; T\in{\cal T}_{sp}(G)) \]
\[ ((p=|X|\,\wedge\, q=|B|\,\wedge\, r=|B|-|B_{T}|)\;\rightarrow\;(r=q-p+1)) \]
By transfer, we obtain the following formula in hypernatural numbers.
\[ (\forall\; p,q,r\in\,^{*}\! \N)\;
(\forall\; ^{*}\! G\in\,^{*}\! {\cal C}_{f})\;
(\forall\; ^{*}\!T\in\,^{*}\! {\cal T}_{sp}( ^{*}\! G)) \]
\[ ((p=| ^{*}\! X|\,\wedge\, q=| ^{*}\! B|\,\wedge\, 
r=| ^{*}\! B|-| ^{*}\! B_{T}|)\;\rightarrow \;(r=q-p+1)) \]

Another standard formula for a connected finite graph having no multiedges
is that $p-1\leq q\leq p(p-1)/2$.  Symbolically, we have
\[ (\forall\; p,q\in\N)\;(\forall\; G\in{\cal C}_{f})\;
((p=|X|\,\wedge\, q=|B|)\;\rightarrow\; (p-1\leq q\leq p(p-1)/2)). \]
So, by transfer, we have 
\[ (\forall\; p,q\in\,^{*}\! \N)\;(\forall\; G\in\,^{*}\! {\cal C}_{f})\;
((p=|^{*}\! X|\,\wedge\, q=|^{*}\! B|)\;\rightarrow\;(p-1\leq q\leq p(p-1)/2)). \]

Still another example of such a lifting concerns the radius $R$ and 
diameter $D$ of a finite connected graph $G$.  It is a fact that
$R\leq D\leq 2R$ \cite[page 37]{b-h}, \cite[pages 20-21]{c-l}. Again, we need to express
ideas symbolically.

Let $A\subset \N$ be such that $|A|\in\N$ (i.e., $A$ is a finite subset of
$\N$).  We use the symbols $\overline{a}=\max A$ and $\underline{a}=\min A$
as abbreviations for the following sentences.
\[ \overline{a}=\max A\; \leftrightarrow\; (\forall\; c\in A)\;(\exists\;\overline{a}\in A)\;
(\overline{a}\geq c) \]
\[ \underline{a}=\min A\; \leftrightarrow\; (\forall\; c\in A)\;(\exists\;\underline{a}\in A)\;
(\underline{a}\leq c) \]
This transfers to
\[ \overline{a}=\max \,^{*}\! A\; \leftrightarrow\; (\forall\; c\in\,^{*}\! A)\;(\exists\;\overline{a}\in\,^{*}\! A)\;
(\overline{a}\geq c) \]
\[ \underline{a}=\min\, ^{*}\! A\; \leftrightarrow\; (\forall\; c\in  ^{*}\! A)\;(\exists\;\underline{a}\in\,^{*}\! A)\;
(\underline{a}\leq c), \]
where now $ ^{*}\! A$ is a hyperfinite subset of $ ^{*}\! \N$ 
and $\overline{a}$ and $\underline{a}$ are hypernatural numbers.
$ ^{*}\! A$ does have a maximum element and a minimum element
so that these definitions are valid \cite[pages 149-150]{go}.

Now, for a given finite connected graph $G\in{\cal C}_{f}$ 
and with $x$ and $y$ being any two vertices in $G$, let
${\cal P}_{x,y}$ denote the set of all paths in $G$ terminating
at $x$ and $y$.  The length $|P_{x,y}|$ of any $P_{x,y}\in{\cal P}_{x,y}$
is the number of edges in $P_{x,y}$.
Then, the distance 
between any two vertices $x$ and $y$ of $G=\{X,B\}$ is the natural number
$\min\{|P_{x,y}|\!: P_{x,y}\in{\cal P}_{x,y}\}$.  So, with 
${\cal D}(G)$ denoting the set of distances in $G$, we may write
\[ d_{x,y}\in{\cal D}(G)\;\leftrightarrow\;(*\exists x,y\in X)
(\forall P_{x,y}\in{\cal P}_{x,y})(d_{x,y}=\min\{|P_{x,y}|\in{\cal P}_{x,y}\}) \]
If $x=y$, we set $d_{x,x}=0$.
By transfer, we have the enlargement $^{*}\!{\cal P}_{x,y}$ of 
${\cal P}_{x,y}$ consisting of all the nonstandard paths terminating
at the nonstandard vertices $x$ and $y$.  We also have
the internal set $^{*}\!{\cal D}(^{*}\!G)$ of 
all hypernatural distances in $^{*}\!G$:
\[ ^{*}\!d_{x,y}\,\in^{*}\!{\cal D}(^{*}\!G)\;\leftrightarrow\;
(\exists x,y\in\,^{*}\!X)(\forall \,^{*}\!P_{x,y}\in\,^{*}\!{\cal P}_{x,y})
(^{*}\!d_{x,y}=\min\{|^{*}\!P_{x,y}|\!:\,^{*}\!P_{x,y}\in\,
^{*}\!{\cal P}_{x,y}\}) \]
Here, too, $^{*}\!d_{x,x}={\bf 0}$.

Furthermore, two nonstandard vertices $x,y\in\,^{*}\!X$ are said to be 
{\em limitedly distant} if $^{*}\!d_{x,y}$ is a limited hypernatural.
This defines an equivalence relation
on $^{*}\!X$; indeed, reflexivity and symmetry are obvious, 
and the triangle inequality follows by transfer from the triangle 
inequality for distances in $G$. So, by analogy with galaxies for hyperreals,
we might refer to the equivalence classes induced by limited distances
as ``nodal galaxies'';  these partition $^{*}\!X$.

Next, we have the eccentricity $e_{x}$ of any vertex of $G=\{X,B\}$
as the maximum of the distances from $x$.  With ${\cal E}(G)$ 
being the set of eccentricities in $G$ and with ${\cal D}_{x}$
being the set of distances from $x$, we have symbolically
\[ e_{x}\in{\cal E}(G)\;\leftrightarrow\;(exists x\in X)
(e_{x}=\max\{d_{x}\!: d_{x}\in{\cal D}_{x}\}). \]
So, for a hyperfinite connected graph $^{*}\!{\cal G}=\{^{*}\!X,\,^{*}\!B\}
\in\,^{*}\!{\cal C}_{f}$, we have by transfer
\[ e_{x}\in{\cal E}(^{*}\!G)\;\leftrightarrow\;(\exists x\in\,^{*}\!X)
(e_{x}=\max\{d_{x}\!: d_{x}\in\,^{*}\!{\cal D}_{x}\}), \]
where now ${\cal E}(^{*}\!G)$ is the set of hypernatural eccentricities 
for nonstandard vertices in 
$^{*}\!G$ and $^{*}\!{\cal D}_{x}$ is the set
of hypernatural distances starting at the nonstandard vertex $x$.

Then, for any $G=\{X,B\}\in {\cal C}_{f}$, the radius $R(G)$ is defined by 
\[ (\forall\; e_{x}\in{\cal E}(G))\;(\exists\; R(G)\in\N)\;
(R(G)=\min\{ e_{x}\!: x\in X\}), \]
which by transfer gives the following definition of the hypernatural
radius $R( ^{*}\! G)\in\,^{*}\! \N$ of any $ ^{*}\! G=\{ ^{*}\! X, \,^{*}\! B\}
\in\,^{*}\! {\cal C}_{f}$:
\[ (\forall\; e_{x}\in\,^{*}\! {\cal E}(G))\;(\exists\; R( ^{*}\! G)\in\,^{*}\! \N)\;
(R( ^{*}\! G)=\min\{ e_{x}\!: x\in\,^{*}\! X\}), \]
Similarly, the diameter $D(G)$ of $G$ is defined by
\[ (\forall\; e_{x}\in{\cal E}(G))\;(\exists\; D(G)\in\N)\;
(D(G)=\max\{ e_{x}\!: x\in X\}), \]
which by transfer gives the hypernatural diameter $D(G)\in\,^{*}\! \N$ of 
$ ^{*}\! G$.
\[ (\forall\; e_{x}\in{\cal E}( ^{*}\! G))\;(\exists\; D( ^{*}\! G)\in\,^{*}\! \N)\;
(D( ^{*}\! G)=\max\{ e_{x}\!: x\in\,^{*}\! X\}), \]
So, we have the following sentence for the standard result:
\[ (\forall\; G\in{\cal C}_{f})\;(R(G)\leq D(G)\leq 2R(G)), \]
which by transfer yields the nonstandard result
\[ (\forall\; G\in\,^{*}\! {\cal C}_{f})\;(R( ^{*}\! G)\leq D( ^{*}\! G)\leq 2R( ^{*}\! G)). \]

9. Eulerian Graphs

A finite trail is defined much as a finite path is defined 
except that the condition that all the vertices be distinct is 
relaxed;  however, edges are still required to be distinct.
Thus, the truth of the following sentence defines a {\em trail} $T$ in a finite 
graph $G=\{X,B\}$, with $T$ having two or more edges.  This time we use 
the notation $b_{m}=\{x_{m},y_{m}\}$ to display the vertices
$x_{m}$ and $y_{m}$ that are incident to $b_{m}$.
\[ (\exists\; k\in\N\setminus\{0\})\;(\exists\; b_{0},b_{1},\ldots,b_{k}\in B)\;
(\forall m\in\{0,\ldots,k-1\})\;(y_{m}=x_{m+1}) \]
That $B$ is a set insures that the edges $b_{0},b_{1},\ldots,b_{k}$ 
are all distinct.  On the other hand, this sentence allows vertices 
to repeat in a trail.

For a {\em closed trail}, we have the truth of the following 
sentence as its definition.
\[ (\exists\; k\in\N\setminus\{0\})\;(\exists\; b_{0},b_{1},\ldots,
b_{k}\in B)\;(\forall m\in\{0,\ldots,k-1\})\;(y_{m}=x_{m+1})\;\wedge\;
(y_{k}=x_{0}) \]
With $Q$ denoting a trail, we denote the set of edges in $Q$
by $B(Q)$.  Also, we let ${\cal Q}(G)$ denote the set of 
closed trails in a given graph $G=\{X,B\}$.

By attaching asterisks as usual, we obtain by transfer the corresponding
sentence for trails in a given nonstandard graph
$ ^{*}\! G=\{ ^{*}\! X,\, ^{*}\! B\}$.  Thus, a {\em nonstandard trail} $ ^{*}\! Q$
is defined by the truth of the following sentence; now, $b_{m}=\{x_{m},y_{m}\}$ 
is a nonstandard edge with the nonstandard vertices $x_{m}$ and $y_{m}$.

\[(\exists\; k\in\,^{*}\! \N\setminus\{0\})\;(\exists\; b_{0},b_{1},\ldots,
b_{k}\in\,^{*}\! B)\;(\forall m\in\{0,\ldots, k-1\})\;(y_{m}=x_{m+1}). \]
A similar expression holds for a {\em nonstandard closed trail} (just append 
$\wedge\;(y_{k}=x_{0})$).
With $ ^{*}\! Q$ denoting a nonstandard trail, we denote
the set of nonstandard edges in $ ^{*}\! Q$ by 
$ ^{*}\! B( ^{*}\! Q)$.  Also, we let $ ^{*}\! {\cal Q}( ^{*}\! G)$ 
denote the set of nonstandard closed trails in a given nonstandard
graph $ ^{*}\! G=\{ ^{*}\! X, ^{*}\! B\}$.

A finite connected graph $G=\{X,B\} \in {\cal C}_{f}$ is called {\em Eulerian} 
if it contains a closed trail that meets every vertex of $X$.  
The degree $d_{x}$ of $x\in X$ is the 
natural number $d_{x}=|\{b\in B\!: x\in b\}|$.  The nonstandard
version of this definition is as follows:
Given $ ^{*}\! G=\{ ^{*}\! X,\, ^{*}\! B\}\in\,^{*}\! {\cal C}_{f}$,
for any $x\in\,^{*}\! X$, the {\em degree} of $x$ is
$d_{x}=|\{b\in\,^{*}\! B\!: x\in b\}|$.  In this case,
$d_{x}$ may be an unlimited hypernatural number 
when $ ^{*}\! G$ is a hyperfinite 
graph.  However, $ ^{*}\! G$ might happen to be a finite 
graph $G\in{\cal C}_{f}$, which from the point of view of an 
ultrapower construction can occur if all the $G_{n}$ for 
$ ^{*}\! G=[G_{n}]$ are the same finite graph $G\in{\cal G}_{f}$;
in this case, $d_{x}$ will be a natural number for all $x\in\,^{*}\! X$.

Let ${\cal E}_{u}$ (resp. $ ^{*}\! {\cal E}_{u}$) 
denote the set of all standard Eulerian graphs
(resp. nonstandard Eulerian graphs).  Then, Eulerian graphs can be
defined by asserting the truth of the following sentence to the
right of $\leftrightarrow$, where as usual $G=\{X,B\}$.
\[ G\in {\cal E}_{u}\; \leftrightarrow\; (\exists\; Q\in{\cal Q}(G))\;
((\forall\; b\in B)\;(b\in B(Q)) \]
By transfer, the truth of the following right-hand side
defines {\em nonstandard Eulerian graphs}.  Now, $^{*}\! G=
\{ ^{*}\!X,\,^{*}\!B\}$.
\[  ^{*}\! G\in\,^{*}\! {\cal E}_{u}\; \leftrightarrow\;(\exists\; ^{*}\! Q\in \,^{*}\! {\cal Q}(G))\;((\forall\; b\in\,^{*}\! B)\;(b\in\,^{*}\! B( ^{*}\! Q))) \]

Now an ancient theorem of Euler asserts that a graph $G$ 
is Eulerian if and only if the degree of every vertex of 
$G=\{X,B\}$ is an even natural number.  Symbolically,
this can be stated as follows.
\[ G\in {\cal E}_{u}\; \leftrightarrow\;(\forall\; x\in X)\;(d_{x}/2\in\N) \]
Transferring this, we get the nonstandard version of this theorem of Euler:
\[  ^{*}\! G\in\,^{*}\! {\cal E}_{u}\; \leftrightarrow\;
(\forall\; x\in  ^{*}\! X)\;(d_{x}/2\in\,^{*}\! \N) \]

10. Hamiltonian Graphs

In this section, it is assumed that each graph $G=\{X,B\}$ is connected 
and finite and
has at least three vertices (i.e., $|X|\geq 3$). A graph $G$
is called {\em Hamiltonian} if it contains a loop that meets every
vertex in the graph.  Let ${\cal L}(G)$ denote the set
of all loops in $G$.  Also, for any loop $L\in{\cal L}(G)$, let $X(L)$
denote the vertex set of $L$.  Then, a 
Hamiltonian graph $G=\{X,B\}\in{\cal C}_{f}$ is also defined by 
the truth of the following sentence to the right of 
$\leftrightarrow$.  The set of Hamiltonian 
graphs will be denoted by $\cal H$,
where ${\cal H}\subset{\cal C}_{f}$.
\[ G\in{\cal H}\; \leftrightarrow\;(\exists\; L\in{\cal L}(G))\;((\forall\; x\in X)\;
(x\in X(L))) \]
By transfer of this sentence, we define a {\em nonstandard Hamiltonian graph}
as follows, where now 
$ ^{*}\! {\cal H}$ is the set of nonstandard hyperfinite Hamiltonian graphs,
${\cal L}( ^{*}\! G)$ is the set of all nonstandard loops in 
$ ^{*}\! G$, $ ^{*}\! X((L)$ is the set of nonstandard vertices
in $L\in{\cal L}( ^{*}\! G)$, and $^{*}\!G=\{^{*}\!X,\,^{*}\!B\}$.
\[  ^{*}\! G\in\,^{*}\! {\cal H}\; \leftrightarrow\;
(\exists\; L\in{\cal L}( ^{*}\! G))\;((\forall\; x
\in\,^{*}\! X)\;(x\in\,^{*}\! X(L))). \]

A simple criterion for a graph $G=\{X,B\}$ to be Hamiltonian is that the 
degree $d_{x}$ of each of its vertices $x$ be no less than one half of $|X|$
\cite[page 134]{b-c}, \cite[page 79]{b-h}.  
Symbolically, this condition is expressed as
follows:
\[ ((\forall\; x\in X)\;(d_{x}\geq |X|/2))\;\rightarrow\;{\cal G}\in{\cal H}. \]
By transfer, we get the following 
criterion for a nonstandard Hamiltonian graph. 
\[ ((\forall\; x\in\,^{*}\! X)\;(d_{x}\geq| ^{*}\! X|/2))\;\rightarrow\; ^{*}\! {\cal G}\in\,^{*}\! {\cal H}. \]

A more general criterion due to Ore asserts that $G=\{X,B\}$
is Hamiltonian if, for every pair of nonadjacent vertices $x$ and $y$, 
$d_{x}+d_{y}\geq|X|$ \cite[page 134]{b-c}, \cite[page 79]{b-h}.
Symbolically, we have
\[ ((\forall\; x,y\in X)\;((\neg(x\diamond y))\rightarrow (d_{x}+d_{y}\geq |X|)))
\;\rightarrow\; G\in{\cal H}, \]
which by transfer becomes
\[ ((\forall\; x,y\in  ^{*}\! X)\;((\neg(x\diamond y))\rightarrow (d_{x}+d_{y}\geq | ^{*}\! X|)))
\;\rightarrow\;  ^{*}\! G\in\,^{*}\! {\cal H}. \]

Still more general is Posa's theorem \cite[page 132]{b-c}, \cite[page 79]{b-h}:
If, for every $j\in\N$ satisfying $1\leq j <|X|/2$, the number of vertices
of degree no larger than $j$ is less than $j$, then 
the graph $G=\{X,B\}$ is Hamiltonian.  The following symbolic 
sentence states this criterion.
\[ ((\forall\; j\in\N)\;(\forall\; x\in X)\;((1\leq j <|X|/2)\rightarrow
(|\{x\in X\!: d_{x}\leq j \}| < j)))\;\rightarrow\; G\in{\cal H} \]
By transfer the following criterion holds for nonstandard
graphs $ ^{*}\! G=\{ ^{*}\! X,\, ^{*}\! B\}$.
\[ ((\forall\; j\in\,^{*}\! \N)\;(\forall\; x\in\,^{*}\! X)\;((1\leq j <| ^{*}\! X|/2)\rightarrow
(|\{x\in\,^{*}\! X\!: d_{x}\leq j \}| < j)))\;\rightarrow\;  ^{*}\! G\in\,^{*}\! {\cal H} \]

11. A Coloring Theorem

A simple graph-coloring theorem that is not restricted to planar graphs
asserts that, if the largest of the degrees for the vertices of a graph 
$G=\{X,B\}$ is $k$, then the graph is $(k+1)$-colorable
\cite[page 82]{wi}.\footnote{There exists a function $f$ that assigns to each
vertex one of $k+1$ colors such that no two adjacent vertices have the 
same color.} To express this symbolically, first let 
$M(X,\N_{k+1})$ denote the set of all functions that map a set $X$ 
into the set $\N_{k+1}$ of those natural numbers $j$ satisfying
$1\leq j\leq k+1$.  Then, the following restates this theorem for the 
given graph $G$.
\[ (\exists\; k\in\N)\;(\forall\; x,y\in X) \]
\[((d_{x}\leq k)\;\rightarrow\; 
((\exists\; f\in M(X,\N_{k+1}))\;(( x\diamond y)
\rightarrow (f(x)\neq f(y))))) \]

To transfer this, we first let $ ^{*}\! M( ^{*}\! X,\N_{k+1})$
be the set of all internal functions mapping the enlargement 
$ ^{*}\! X$ into $\N_{k+1}$.  Then, this theorem is transferred
to nonstandard graphs simply by appending asterisks, as usual:
\[ (\exists\; k\in\N)\;(\forall\; x,y\in\,^{*}\! X) \]
\[ ((d_{x}\leq k)\;\rightarrow\; 
((\exists\;  ^{*}\! f\in\,^{*}\! M (^{*}\! X,\N_{k+1}))\;(( x\diamond y)
\rightarrow ( ^{*}\! f(x)\neq\, ^{*}\! f(y))))) \]
Note here that the assumption of a natural-number 
bound $k$ on the degrees of all 
the nonstandard vertices has been maintained.  This conforms to the 
fact that the enlargement of the finite set $\N_{k+1}$ is 
$\N_{k+1}$.  As a consequence, the conclusion remains strong.  

On the other hand, we could generalize this transferred theorem
as follows: In terms 
of an ultrapower construction,
we could replace $\N_{k+1}$ by an internal set $ ^{*}\! \N_{k+1}$
obtained from a sequence $\langle \N_{k_{n}+1}\!: n\in \N\rangle$
of finite sets $\N_{k_{n}+1}$, one set for each $G_{n}$
with regard to $ ^{*}\! G=[G_{n}]$ and with $k_{n}$ being the maximum
vertex degree in $G_{n}$.  
But then, our conclusion would be
weakened to a coloring with a hypernatural number $ ^{*}\! k
=[k_{n}]$ of colors.

12. A Final Comment

Undoubtedly, other standard results for graphs can be 
lifted in this way to nonstandard settings.

\end{document}